\documentclass[a4paper,12pt]{article}
\usepackage[french]{babel}
\usepackage{latexsym,amssymb}
\usepackage{enumerate}
\usepackage{amsmath}
\usepackage{pb-diagram}
\newtheorem{thm}{Theorem}[section]

\newtheorem{pro}[thm]{Proposition}

\newcommand{\HH}{\mathbb{H}}
\newcommand{\C}{\mathbb{C}}

\newcommand{\R}{\mathbb{R}}
\newcommand{\T}{\mathbb{T}}
\newcommand{\s}{\mathbb{S}}

\def \be{\begin{eqnarray*}}
\def \ee{\end{eqnarray*}}
\def \ben{\begin{enumerate}}
\def \een{\end{enumerate}}
\def \beit{\begin{itemize}}
\def \eeit{\end{itemize}}
\def \bui#1#2{\mathrel{\mathop{\kern 0pt#1}\limits^{#2}}}
\def \buil#1#2{\mathrel{\mathop{\kern 0pt#1}\limits_{#2}}}
\def \bflr{\begin{flushright}}
\def \eflr{\end{flushright}}

\begin{document}
\title{Remarques sur les spineurs de Killing transversaux}
\author{{\small{Nicolas Ginoux}}
\footnote{Institut f\"ur Mathematik - Geometrie,
Universit\"at Potsdam,
Am Neuen Palais 10,
D-14469 Potsdam,
E-mail: ginoux@math.uni-potsdam.de}, {\small Georges Habib}\footnote{Max-Planck-Institut f\"ur Mathematik in den Naturwissenschaften,
Inselstr. 22,
D-04103 Leipzig,
E-mail: Georges.Habib@mis.mpg.de}}
\date{}
\maketitle

\noindent {\bf R\'esum\'e} Nous montrons des conditions n\'ecessaires d'existence de spineurs de Killing transversaux sur une vari\'et\'e spinorielle munie d'un flot riemannien. \\ 

\begin{center}{\large Remarks on transversal Killing spinors}
\end{center}

\noindent {\bf Abstract} We show necessary conditions for the existence of transversal Killing spinors on a spin manifold endowed with a Riemannian flow.\\ 

\noindent Soit $(M,g,\mathcal{F})$ une vari\'et\'e riemannienne munie d'un feuilletage riemannien $\mathcal{F}$. Notons par $Q$ le fibr\'e normal du feuilletage et $\nabla$ la connexion de Levi-Civita transversale. Dans \cite{Hab}, nous d\'efinissons la notion de spineurs de Killing transversaux, i.e. les champs de spineurs basiques qui satisfont pour tout $Z\in \Gamma(Q)$ l'\'equation $\nabla_{Z}\psi=\beta Z\cdot\psi$    
o\`u $\beta$ est un nombre complexe et $``\cdot"$ d\'esigne la mutiplication de Clifford du fibr\'e des spineurs de $Q$, suppos\'e spinoriel. Ces spineurs se manifestent naturellement dans l'\'etude du spectre de l'op\'erateur de Dirac basique. Dans \cite{GinHab}, nous \'etudions les conditions d'int\'egrabilit\'e de ces spineurs sur les flots riemanniens (feuilletages de dimension $1$) et nous donnons une classification partielle sur les vari\'et\'es de Sasaki et les flots de dimension $3$. \\\\
Dans cette note, nous montrons d'une part la non-existence de spineurs de Killing transversaux imaginaires sur une vari\'et\'e compacte, d'autre part nous \'etendons la classification en dimension $3$ sans hypoth\`eses de minimalit\'e ni de courbure.

\noindent Les auteurs aimeraient remercier le SFB 647 ainsi que l'institut de Max-Planck pour les Math\'ematiques dans les Sciences de Leipzig. C'est aussi un grand plaisir de remercier Oussama Hijazi et Christian B\"ar pour leur soutien.  

\section{Pr\'eliminaires} 

\noindent Dans cette partie, on va rappeler quelques \'el\'ements de base des feuilletages spinoriels, on se ref\`ere \`a \cite{Car84, Ton, Kac90, Hab}. Soit $(M,g,\mathcal{F})$ une vari\'et\'e riemannienne de dimension $n$ munie d'un feuilletage $\mathcal{F}$ donn\'e par un sous-fibr\'e $L\subset TM$. Le feuilletage $\mathcal{F}$ est dit riemannien si la m\'etrique $g$ est quasi-fibr\'ee, i.e. la m\'etrique induite de $g$ sur le fibr\'e normal $Q:=L^\perp,$ de rang $q$, satisfait la condition d'invariance d'holonomie. Cette condition est caract\'eris\'ee par l'existence d'une connexion m\'etrique \`a torsion nulle $\nabla$ sur $Q$. La courbure transversale $R^\nabla$ est une forme basique dans le sens qu'elle s'annulle le long des feuilles. On peut ainsi associer \`a $\nabla$ la courbure de Ricci transversale ${\rm Ric}^\nabla$ et sectionelle transversale $K^\nabla$. Dans toute la suite, on va consid\'erer des feuilletages riemanniens.\\
Supposons maintenant que $Q$ porte une structure spinorielle (en tant qu'un fibr\'e vectoriel feuillet\'e) et notons par $\Sigma Q$ son fibr\'e des spineurs. L'op\'erateur de Dirac basic $D_{b}$ agit sur les sections basiques du fibr\'e des spineurs $\Sigma Q$ pour tout $\psi$ par 
\begin{equation}
D_{b}\psi=\sum_{i=1}^q e_i\cdot \nabla_{e_i}\psi-\frac{1}{2}\kappa\cdot\psi,
\label{eq:bas}
\end{equation}
o\`u $\{e_i\}_{i=1,\cdots,q}$ est un rep\`ere local de $Q$ et $\kappa$ est la courbure moyenne des feuilles. Ici, $``\cdot"$ d\'esigne l'action de $Q$ sur le fibr\'e des spineurs $\Sigma Q$ par la multiplication de Clifford. Cet op\'erateur est transversalement elliptique et auto-adjoint si $M$ est compacte et donc a un spectre discret \cite{Kac90}. 

\section{R\'esultats}            	
\noindent Soient $(M,g,\mathcal{F})$ une vari\'et\'e riemannienne et $\mathcal{F}$ un feuilletage spinoriel. On rappelle, que si $M$ est compacte, la courbure moyenne $\kappa$ est une $1$-forme basique et harmonique \cite{MM96}. Supposons que le fibr\'e normal admet un spineur Killing transversal pour $\beta\in \C$. En utilisant la formule de Ricci transversale, on montre que le feuilletage est transversalement Einstein et on a ${\rm Ric}^\nabla Z=4(q-1)\beta^2 Z$ pour tout $Z\in \Gamma(Q)$ \cite{Hab}. Donc une premi\`ere cons\'equence est que $\beta$ est r\'eel ou imaginaire pur.  

\begin{pro}\label{Killima} Soit $(M,g,\mathcal{F})$ une vari\'et\'e riemannienne compacte munie d'un feuilletage spinoriel $\mathcal{F}$ et admettant un spineur de Killing transversal non nul pour la constante $\beta$. Alors $\beta\in \R$ et, si de plus $\beta\neq 0,$ le feuilletage est minimal.
\end{pro}
{\bf Preuve.} Supposons d'abord, par l'absurde, que $\beta\in i\R^*$. En appliquant l'op\'erateur de Dirac basique \`a $\psi$, on obtient par \eqref{eq:bas} que $D_b\psi=-q\beta\psi-\frac{1}{2}\kappa\cdot\psi$. L'int\'egrale sur $M$ du produit hermitien de cette identit\'e avec $\psi$ donne, apr\`es avoir identifi\'e les parties imaginaires pures $\beta=-\frac{1}{2q} \frac{\int_M(\kappa\cdot\psi,\psi)v_g}{\int_M|\psi|^2v_g}$.
Maintenant on va calculer le terme $D_b^2\psi$ \`a l'aide de la formule de Schr\"odinger-Lichnerowicz transversal. En effet, on a  
$\nabla^*\nabla\psi=q\beta^2\psi+\beta\kappa\cdot\psi$
et donc on trouve d'une part 
$D_b^2\psi=(q^2\beta^2+\frac{1}{4}|\kappa|^2)\psi+\beta\kappa\cdot\psi$, et d'autre part
$|D_b\psi|^2=q^2|\beta|^2|\psi|^2+q(\kappa\cdot\psi,\beta\psi)+\frac{1}{4}|\kappa|^2|\psi|^2$.
En comparant l'int\'egrale sur $M$ de cette derni\`ere identit\'e avec $\int_M(D_b^2\psi,\psi)v_g$, on en d\'eduit que
$$\int_M[(q^2\beta^2+\frac{1}{4}|\kappa|^2)|\psi|^2+\beta(\kappa\cdot\psi,\psi)]v_g=\int_M[q^2|\beta|^2|\psi|^2+q(\kappa\cdot\psi,\beta\psi)+\frac{1}{4}|\kappa|^2|\psi|^2] v_g.$$ 
Finalement, en \'ecrivant $\beta=iy=\frac{-i}{2q}\frac{\int_M{\rm Im}(\kappa\cdot\psi,\psi)v_g}{\int_M|\psi|^2v_g}$ et en la rempla\c cant dans l'\'egalit\'e ci-dessus on trouve  
que $qy^2=0$ et on arrive \`a une contradiction. Pour $\beta\in \R^*$ la positivit\'e de la courbure de Ricci transversale implique l'annulation de la premi\`ere classe de cohomologie basique $H^1(M/\mathcal{F})$ de $\mathcal{F}$ \cite{Heb} et le r\'esultat se d\'eduit par \cite[Cor. 4.8]{Kac97}.
\hfill$\square$\\


\noindent On suppose d\'esormais que $M$ est de dimension $3$ et est munie d'un flot riemannien, i.e. d'un feuilletage de dimension $1$ d\'efini par un champ de vecteurs unitaires $\xi$.

\begin{pro} \label{classification} Soit $(M^3,g,\xi)$ une vari\'et\'e riemannienne compl\`ete munie d'un flot riemannien. Supposons que $M$ admette un spineur de Killing transversal non nul associ\'e \`a la constante $\beta$, alors on a 
\begin{enumerate}
\item Cas o\`u $\beta=0$: Le rev\^etement universel de $M$ est diff\'eomorphe \`a $\R^3$. Si de plus $M$ est compacte, alors c'est le tore plat $\T^3$ ou une fibration de Seifert \`a fibres $\s^1$ avec $\chi(M/\s^1)=0$, et si en outre $M/\s^1$ est spin, alors $M$ est un fibr\'e en cercles sur $\mathbb{T}^2$.  
\item Cas o\`u $\beta\neq 0$: Le rev\^etement universel de $M$ est une fibration sur $\s^2$ (cas o\`u $\beta\in \R$) ou une fibration sur $\HH^2$ (cas o\`u $\beta\in i\R$). Si de plus $M$ est compacte, alors elle est soit le produit riemannien $\s^1\times\s^2$, soit une fibration de Seifert \`a fibres $\s^1$ avec $\chi(M/\s^1)>0$. Si de plus $M/\s^1$ est spin, alors $M$ est un fibr\'e en cercles sur $\s^2$. 
\end{enumerate}
\end{pro} 
{\bf Preuve.} L'existence d'un spineur de Killing transversal associ\'e \`a $\beta$ implique que la courbure sectionnelle transversale est \'egale \`a $4\beta^2$. Quitte \`a effectuer une homoth\'etie sur la m\'etrique transversale on peut supposer que $\beta\in\{\pm\frac{1}{2},0,\pm\frac{i}{2}\}$. Lorsque $\beta=0$ le fibr\'e normal de $\mathcal{F}$ est plat, par cons\'equent le rev\^etement universel $\widetilde{M}$  de $M$ est diff\'eomorphe \`a $\R^3$ \cite[Cor. 3.2]{Blu84}. Pour $\beta\neq 0$, la courbure transversale $R^\nabla$ v\'erifiant $\nabla R^\nabla=0$, on d\'eduit de \cite{Blu84} que $\widetilde M$ est une fibration riemannienne sur $N=\s^2$ (lorsque $\beta=\pm\frac{1}{2}$) ou sur $N=\HH^2$ (lorsque $\beta=\pm\frac{i}{2}$).
Lorsque $M$ est compacte une version transversale de la formule de Gauss-Bonnet \cite[Prop. 4]{NR} implique que $\chi(M/\s^1)$ est proportionnel \`a la courbure sectionnelle transversale totale sur $M$; par la classification de Y. Carri\`ere \cite{Car84} on en d\'eduit les diff\'erents types de g\'eom\'etrie de $M$. En particulier le quotient $M/\s^1$ est toujours une orbisurface. Si ce quotient est en outre suppos\'e spin, alors c'est une orbisurface orientable dont les singularit\'es - en nombre fini - sont de type $\R^2/\mathbb{Z}_k$, o\`u $k$ est un entier impair \cite{BGR}. D'autre part le spineur de Killing transversal de $M$ descend dans ce cas sur $M/\s^1$. Or la condition d'\'equivariance de ce spineur en chaque point singulier se traduit par la dichotomie suivante: soit cette singularit\'e n'en est pas une (i.e., $k=1$), soit elle est un z\'ero du spineur. Mais un spineur de Killing non trivial n'a pas de z\'ero, par cons\'equent $M/\s^1$ est soit $\T^2$ soit $\s^2$.\hfill$\square$

\end{document}